\documentclass[11pt]{amsart}
\usepackage{amssymb,latexsym,comment,url,amsmath}

\usepackage[T1]{fontenc}
\usepackage{tikz}
\usepackage{rotating,graphicx}
\usepackage{currvita}
\usepackage{diagbox} 
\usepackage{xcolor}
\usepackage{subcaption}
\usepackage{multirow}
\usepackage{array}
\usepackage{longtable}
\usepackage{ltablex}
\usepackage{float}
\usepackage{lmodern,babel,adjustbox,booktabs,multirow}
\usepackage[T1]{fontenc}
\usepackage[utf8]{inputenc}
\usepackage{ upgreek }
\usepackage{caption}
\usepackage{placeins}
\restylefloat{table}
\usepackage{blindtext}
\usepackage{hyperref}
\usepackage{verbatim}
\usepackage{verbatimbox}

\newcommand{\Q}{{\mathbb Q}}

\newcommand{\Z}{{\mathbb Z}}

\newcommand{\Magma}{{\sf Magma }}

\newenvironment{Proof}{\par\noindent{\sc Proof:}}%
                      {\hspace*{\fill}\nobreak$\Box$\par\medskip}
                       {\hspace*{\fill}\nobreak$\Box$\par\medskip}

\newtheorem{Proposition}{Proposition}[section]
\newtheorem{Theorem}[Proposition]{Theorem}

\makeglossary
\newtheorem*{theorem*}{Theorem}

\theoremstyle{definition}

 \newtheorem{Example}[Proposition]{Example}

\renewcommand{\baselinestretch}{1.1}

\begin{document}

\title[Families of twists of tuples of hyperelliptic curves]%
{Families of twists of tuples of hyperelliptic curves}

\author[B. M. Amir]%
{Beyza Mevl\"{u}de Amir}
\address{Faculty of Engineering and Natural Sciences, Sabanc{\i} University, Tuzla, \.{I}stanbul, 34956 Turkey}
\email{beyza.cepni@sabanciuniv.edu}
\author[M. Sadek]%
{Mohammad~Sadek}
\address{Faculty of Engineering and Natural Sciences, Sabanc{\i} University, Tuzla, \.{I}stanbul, 34956 Turkey}
\email{mohammad.sadek@sabanciuniv.edu}
\author[N. El-Sissi]%
{Nermine El-Sissi}
\address{Bah\c{c}e\c{s}ehir University, Faculty of Engineering
and Natural Sciences, Istanbul, Turkey}
\email{nermineahmed.elsissi@eng.bau.edu.tr }

\begin{abstract}
Let $f \in \Q[x]$ be a square-free polynomial of degree at least $3$,
$m_i$, $i=1,2,3$, odd positive integers, and $a_i$, $i=1,2,3$, non-zero rational numbers. We show the existence of a rational function $D\in\Q(v_1,v_2,v_3,v_4)$ such that the Jacobian of the quadratic twist of $y^2=f(x)$ and the Jacobian of the $m_i$-twist, respectively $2m_i$-twist, of $y^2=x^{m_i}+a_i^2$, $i=1,2,3$, by $D$ are all of positive Mordell-Weil ranks. As an application, we present families of hyperelliptic curves with large Mordell-Weil rank.
\end{abstract}
\maketitle

\let\thefootnote\relax

\section{Introduction}
Given an abelian variety $A$ defined over a number field $K$, the Mordell-Weil Theorem asserts that the abelian group of rational points $A(K)$ is finitely generated. This implies that $A(K)$ is isomorphic to $\mathbb{T}\times\Z^r$, where $\mathbb T$ is the torsion subgroup of $A(K)$, and $r$ is a non-negative integer, called the Mordell-Weil rank of $A(K)$. In this article, we are interested in the following question. Fix an integer $n\ge2$. Given (hyper)elliptic curves $C_i$ defined over $\Q$ of genus $g_i$, and tuples of positive integers $r_i$, $i=1,\cdots,n$, is there a twist of $C_i$ such that the Mordell-Weil rank of the Jacobian of this twist is at least $r_i$, $i=1,\cdots,n$? 

The problem was initially discussed in \cite{Kuwata} where it was proved that for a pair of elliptic curves $E_1$ and $E_2$ over $\Q$ whose $j$-invariants are not simultaneously $0$ or $1728$, there exist infinitely many square-free rational numbers $d$ such that the quadratic twists of $E_1$ and $E_2$ by $d$ are both of positive Mordell-Weil rank. In \cite{Alaa}, it was proved that there exist families of pairs of elliptic curves $E_1$ and $E_2$ and infinitely many square-free rationals $d$ such that the quadratic twists of $E_1$ and $E_2$ by $d$ are both of
Mordell-Weil rank at least $2$. The positivity of the Mordell-Weil rank of families of twists of triples and quadruples of elliptic curves with $j$-invariants $0$ or $1728$ was established in \cite{Ulas2010G,Ulas2013}. This motivates the investigation of the positivity of Mordell-Weil rank of families of twists of higher dimensional abelian varieities.

Given non-zero rational numbers $a, b, c$, it was proved that there exists
a polynomial $d(t) \in \Q[t]$ such that the Jacobians of the curves given by $y^2 = x^n + ad(t)$, $y^2 = x^n + bd(t)$, $y^2=x^n+cd(t)$, $n\ge3$, both have positive Mordell-Weil rank over $\Q(t)$, see \cite{Ulas2011}. The latter result was extended to four such hyperelliptic curves when $n$ is odd. 
If $f \in \Q[x]$ is a square-free polynomial of degree at least $3$, it was proved in \cite{Ulas2014} that
 there exists a function $D \in
\Q(u, v, w)$ such that the Jacobians of the curves
 $Dy^2 = f(x)$, $y^2 = Dx^m + b$, and $C$ all
have positive Mordell-Weil ranks over $\Q(u, v, w)$ where $C$ is either $y^2=Dx^m+c$ or $y^2=x^m+cD$, and $m\ge 3$ is odd. The question was investigated for other families of quadruples of hyperelliptic curves in \cite{SadekMesut}.

In this work, we extend the latter results to include several families of quadruples of twisted hyperelliptic curves. Namely, we consider the following sets of hyperelliptic curves 
\begin{eqnarray*}
Dy^2=f(x),\quad y^2=Dx^{m_1}+a,\quad y^2=Dx^{m_2}+b,\quad y^2=Dx^{m_3}+c,
\end{eqnarray*}
or
\begin{eqnarray*}
Dy^2=f(x),\quad y^2=x^{m_1}+aD,\quad y^2=x^{m_2}+bD,\quad y^2=x^{m_3}+cD,
\end{eqnarray*}
or
\begin{eqnarray*}
 y^2=Dx^{m_1}+a,\quad y^2=Dx^{m_2}+b,\quad y^2=Dx^{m_3}+c,\quad y^2=Dx^{m_4}+d,
\end{eqnarray*}
where $f(x)\in\Q[x]$ is a square-free polynomial of degree at least $3$, and $m_i\ge 3$ is an odd integer, $i=1,2,3,4$. We prove that if $a$, $b$, $c$ are non-zero rational squares, then $D$ can be chosen to be a rational function in $\Q(u,v_1,v_2,v_3)$ such that the Mordell-Weil ranks of all four hyperelliptic curves over $\Q(u,v_1,v_2,v_3)$ are positive. The proof depends on associating a system of Diophantine equations to these quadruples of hyperelliptic curves. Finding such a rational $D$ corresponds to finding a rational solution to this system of equations.  

As a byproduct, fixing an odd integer $m\ge 3$, we show that there exists a rational $D\in\Q(u,v_1,v_2,v_3)$ such that the Mordell-Weil rank of the Jacobian of $Dy^2=f(x)$ is positive and the Mordell-Weil rank of the Jacobian of $C$  is at least $3$; where $C$ is either $y^2=Dx^m+a^2$ if $\deg f$ is at least $3$; or $y^2=x^m+a^2D$ if $f(x)$ is of degree $3$ or $4$ .

\subsection*{Acknowledgment}
The authors would like to thank the anonymous referee for many useful suggestions and comments  that improved the manuscript. All the calculations in this work were performed using \textbf{Magma}, \cite{Magma}. 
  M. Sadek is supported by The Scientific and Technological Research Council of Turkey, T\"{U}B\.{I}TAK, research grant ARDEB 1001/122F312.

The authors declare no conflicts of interest.

\section{A quadratic twist and three higher degree twists}
Let $f(x)\in\Q[x]$ be a square-free polynomial of degree at least $3$. Fix odd integers $m_1,m_2,m_3\ge 3$. We are investigating the existence of non-square rational numbers $D$ such that the Jacobians of the (hyperelliptic) curves 
\begin{eqnarray*}
    C:Dy^2=f(x),\quad C_1:y^2=Dx^{m_1}+a,\quad C_2:y^2=Dx^{m_2}+b,\quad C_3:y^2=Dx^{m_3}+c
\end{eqnarray*}
are of positive Mordell-Weil rank over $\Q$. 

We set $M=\operatorname{lcm}(m_1,m_2,m_3)$, and $M_{i}=M/m_i$, $i=1,2,3$. The rational number $D$ can be found by solving the following system of equations
\begin{eqnarray}
\label{eq1}
\frac{y_1^2-a}{x_1^{m_1}}=\frac{y_2^2-b}{x_2^{m_2}}=\frac{y_3^2-c}{x_3^{m_3}}=\frac{f(x_4)}{y_4^2}.
\end{eqnarray}
More precisely, we are looking for solutions $x_i,y_i$, $i=1,2,3,4$, for the system (\ref{eq1}) where $$x_1x_2x_3(y_1^2 -a)(y_2^2-b)(y_3^2-c)y_4 f(x_4) \ne0.$$
In what follows, we suggest how one can obtain a family of parametric solutions to (\ref{eq1}). We set 
\begin{eqnarray*}
x_1=\frac{1}{v_1^{2}T^{M_1}},\, x_2=\frac{1}{v_2^2T^{M_2}},\, x_3=\frac{1}{v_3^2T^{M_3}},\, x_4 =u,\,
y_1=p,\, y_2=q,\,y_3=r,\,y_4=\frac{1}{T^{(M-1)/2}}\end{eqnarray*} 
where $u,v_1,v_2,v_3$ are rational parameters. We observe that 
\begin{eqnarray*}
T&=&\frac{f(u)}{v_3^{2m_3}(r^2-c)},\\
v_1^{2m_1}T^M(p^2-a)&=&v_2^{2m_2}T^M(q^2-b)=v_3^{2m_3}T^M(r^2-c).
\end{eqnarray*}

In other words, the system (\ref{eq1}) can be simplified to 

\begin{eqnarray*}
T=\frac{f(u)}{v_3^{2m_3}(r^2-c)},\quad
v_1^{2m_1}(p^2-a)=v_2^{2m_2}(q^2-b)=v_3^{2m_3}(r^2-c).
\end{eqnarray*}

Geometrically, the second equation defines an intersection $\mathcal{C}_{a,b,c}$ of two quadratic surfaces in $\mathbb{P}^3$ over $\Q(v_1,v_2,v_3)$. Since $[p:q:r:s]=[\pm v_2^{m_2}v_3^{m_3}:\pm v_1^{m_1}v_3^{m_3}:\pm v_1^{m_1}v_2^{m_2}:0]$ are rational points on $\mathcal{C}_{a,b,c}$, it follows that $\mathcal{C}_{a,b,c}$ is an elliptic curve over 
$\Q(v_1,v_2,v_3)$. One may check that the latter four projective rational points form a subgroup isomorphic to $\Z/2\Z\times\Z/2\Z$. 

We now assume, further, that $a$, $b$ and $c$ are non-zero rational squares. After a transformation, the curve $\mathcal C_{a,b,c}$ is defined by $$x^2-T_1^2=y^2-T_2^2=z^2-T_3^2,$$ where $T_1=av_1^{m_1}$, $T_2=bv_2^{m_2}$, and  $T_3=cv_3^{m_3}$.

\begin{Proposition}
\label{prop1}
Let $C_{T_1,T_2,T_3}$ be the elliptic curve defined by 
$$x^2-T_1^2= y^2-T_2^2=z^2-T_3^2$$ 
over $\Q(T_1,T_2,T_3)$. The curve $C_{T_1,T_2,T_3}$ has positive Mordell-Weil rank over $\Q(T_1,T_2,T_3)$. In particular, except for a thin set of triples $(t_1,t_2,t_3)\in\Q\times\Q\times\Q$, the curve $C_{t_1,t_2,t_3}$ has positive Mordell-Weil rank over $\Q$.
\end{Proposition}
\begin{Proof}
    The point $[T_1:T_2:T_3:1]\in C_{T_1,T_2,T_3}(\Q(T_1,T_2,T_3))$. In addition, it specializes to a point of infinite order on $C_{1,2,3}$ when $[1:1:1:0]$ is regarded as the identity element, \Magma \cite{Magma}. For the convenience of the reader, we include the code as it will be used later. 

\begin{verbatim}
K<T1,T2,T3>:=FunctionField(Rationals(),3);
P<x,y,z,w>:=ProjectiveSpace(K,3);
C:=Curve(P, [(x^2-T1^2*w^2)-(y^2-T2^2*w^2),(y^2-T2^2*w^2)-(z^2-T3^2*w^2)]);
P:=C![1,1,1,0];
Q:=C![T1,T2,T3];
E,phi:=EllipticCurve(C,P);
Q1:=phi(Q);
\end{verbatim}
    
    Now the statement follows from Silverman's specialization Theorem, \cite[\S 20, Theorem 20.1]{Silverman}.
\end{Proof}
One remarks that when the point $[w_1:w_2:w_3:1]$ is a point of infinite order in $\mathcal C_{a,b,c}$, then any of its multiples can be used to generate a solution to the system (\ref{eq1}). 

 Doubling the point $P:=[a:b:c:1]\in \mathcal{C}_{a^2,b^2,c^2}(\Q(v_1,v_2,v_3))$ yields that $2P=[w_1:w_2:w_3:1]$ is given by 
{\footnotesize\begin{eqnarray*} \Big[ \frac{a^2 b^2 v_1^{2m_1} v_2^{2m_2} + a^2 c^2 v_1^{2m_1} v_3^{2m_3} -
    b^2 c^2 v_2^{2m_2} v_3^{2m_3}}{2a b c v_1^{2m_1} v_2^{m_2} v_3^{m_3}}&:&\frac{a^2 b^2 v_1^{2m_1} v_2^{2m_2} -
    a^2 c^2 v_1^{2m_1} v_3^{2m_3} + b^2 c^2 v_2^{2m_2} v_3^{2m_3}}{2a b c v_1^{m_1} v_2^{2m_2} v_3^{m_3}}\\
    &:& \frac{-a^2 b^2 v_1^{2m_1} v_2^{2m_2} + a^2 c^2 v_1^{2m_1} v_3^{2m_3} +
    b^2 c^2 v_2^{2m_2} v_3^{2m_3}}{2 a b c v_1^{m_1} v_2^{m_2} v_3^{2m_3}}:1 \Big].\end{eqnarray*}}

We now set 
 \begin{eqnarray*}
     T=\frac{f(u)}{v_3^{2m_3}(w_3^2-c^2)},\qquad
     D=(w_1^2-a^2)v_1^{2m_1}T^M.
 \end{eqnarray*} 
 We consider the following twists of hyperelliptic curves 
 $$C: D y^2=f(x),\quad C_1:y^2=Dx^{m_1}+a^2,\quad C_2:y^2=Dx^{m_2}+b^2,\quad C_3:y^2=Dx^{m_3}+c^2.$$
 Now, the points $\displaystyle P_i:=(x_i,y_i)=\left(\frac{1}{v_i^2T^{M_i}},w_i\right)$ are $\Q(v_1,v_2,v_3)$-rational points in $C_i$, $i=1,2,3$, whereas $P:=(x_4,y_4)=\displaystyle\left(u,\frac{1}{T^{(M-1)/2}}\right)$ is a $\Q(u,v_1,v_2,v_3)$-rational point in $C$.

In \cite{Ulas2010G,Ulas2013}, explicit constructions of families of triples and quadruples of elliptic curves with $j$-invariants $0$ or $1728$ were given such that there is a rational function $D$ for which the twists of the curves in each tuple by $D$ are of positive Mordell-Weil rank. This was extended to pairs and triples of Jacobians of hyperelliptic curves of the same genus in \cite{Ulas2011,Ulas2014}. Following the argument above, we obtain the following theorem which can be considered as an extension of the aforementioned results to quadruples of Jacobians of hyperelliptic curves that are not necessarily of the same genus.
 \begin{Theorem}
 \label{thm1}
     Let $f \in \Q[x]$ be a square-free polynomial of degree at least $3$. Let $m_1,m_2,m_3\ge 3$ be odd integers. Consider the hyperelliptic curves 
     $$E:y^2=f(x),\qquad E_1:y^2=x^{m_1}+a^2,\qquad E_2: y^2=x^{m_2}+b^2,\qquad E_3:y^2=x^{m_3}+c^2,$$ where $a,b,c$ are non-zero rational numbers.
 Then there exists a rational function $D_{2,m_1,m_2,m_3} \in \Q(u, v_1,v_2,v_3)$ such that the
Jacobian of the quadratic twist of the curve $E$ by by $D_{2,m_1,m_2,m_3}$ and the Jacobians of the $m_i$-twists of the curves $E_i$, $i=1,2,3$,
 by $D_{2,m_1,m_2,m_3}$ have positive Mordell-Weil rank over $\Q(u, v_1, v_2,v_3)$.
 \end{Theorem}
\begin{Proof}
Let $C$ and $C_i$, $i=1,2,3$, be the twists of $E$ and $E_i$, $i=1,2,3$, described above. Let $P_i$, $i=1,2,3$, and $P$ be the rational points above. We set $J$ and $J_i$ to be the Jacobians of $C$ and $C_i$, $i=1,2,3$, respectively. We consider the rational divisors $D=(P)-(\infty)$ on $C$; and $D_i= (P_i)-(\infty)$, $i=1,2,3$, on $C_i$. Since $C$
is a non-constant quadratic
twist of a constant curve $E$, it follows that the divisor $D$ defines a rational point of infinite
order in $J(\Q(u, v_1, v_2,v_3))$.

Due to \cite[Proposition 2.1]{Ulas2011}, the divisor $D_i$ corresponds to a rational point of infinite order in $J_i(\Q(u, v_1, v_2,v_3))$ as $P_i$, $i=1,2,3$, has non-constant coordinates with a non-zero $y$-coordinate. 
\end{Proof}
\begin{Example}
The curves $E: D y^2=x^5+x+1$, and $E_i: y^2=D x^{4i+1}+i^2$, $i=1,2,3$, are of positive Mordell-Weil rank over $\Q$, where 
\begin{eqnarray*}
D:&=&D(u,v_1,v_2,v_3)=\frac{(u^5+u+1)^{585}(w_1^2-1)v_1^{10}}{v_3^{15210}(w_3^2-9)^{585}}\textrm{ where} \\
w_1&=&\frac{3 v_3^{13}}{4 v_2^9} + v_2^9\left( \frac{1}{3 v_3^{13}} - \frac{3 v_3^{13}}{v_1^{10}}\right),\\
w_2&=&\frac{3 v_3^{13}}{v_1^5} + v_1^5\left( \frac{1}{3 v_3^{13}} - \frac{3 v_3^{13}}{4v_2^{18}}\right),\\
w_3&=&\frac{3 v_2^{9}}{v_1^5} + v_1^5\left( \frac{3}{4 v_2^{9}} - \frac{v_2^{9}}{3v_3^{26}}\right),\quad u,v_1,v_2,v_3\in\Q\setminus\{0\}.
\end{eqnarray*}
Due to Silverman's Specialization Theorem of abelian varieties, one knows that for all but a thin set of values of $u,v_1,v_2,v_3$ in $\Q$, the points $\left(u,\frac{v_3^{7592}(w_3^2-9)^{292}}{(u^5+u+1)^{292}}\right)-(\infty)$ on the Jacobian of $E$, and the point $(P_i)-(\infty)$, $i=1,2,3$, on the Jacobian of $E_i$, where $$P_i=\left(\frac{v_3^{26M_i}(w_3^2-9)^{M_i}}{v_i^2 (u^5+u+1)^{M_i}},w_i\right), \,M_1=117,\,M_2=65,\,M_3=45,$$
are points of infinite order.
\end{Example}

In Corollary 3.2 and Corollary 4.2 of \cite{Ulas2014}, the authors established the existence of pairs of hyperelliptic curves $C_1$ and $C_2$ together with a rational function $D$ such that the Mordell-Weil rank of the Jacobian of a twist of $C_1$ by $D$ is positive whereas the Mordell-Weil rank of the Jacobian of a twist of $C_2$ by $D$ is at least $2$. In the following theorem, we present such examples of pairs of hyperelliptic curves where the Mordell-Weil rank of the Jacobian of a twist of $C_2$ by $D$ is at least $3$.  
\begin{Theorem}
 \label{thm2}
     Let $f \in \Q[x]$ be a square-free polynomial of degree at least $3$. Let $m\ge 3$ be an odd integer. Consider the hyperelliptic curves 
     $$E:y^2=f(x),\qquad E':y^2=x^{m}+a^2,\quad a\in\Q\setminus\{0\}.$$ 
 Then there exists a rational function $D_{2,m} \in \Q(u, v_1,v_2,v_3)$ such that the
Jacobian of the quadratic twist of the curve $E$ by $D_{2,m}$ is of positive Mordell-Weil rank and the Mordell-Weil rank of the Jacobian of the $m_i$-twist of the curve $E'$
 by $D_{2,m}$ over $\Q(u, v_1, v_2,v_3)$ is at least $3$.
 \end{Theorem}
 \begin{Proof}
In Theorem \ref{thm1}, we set $a=b=c$ and $m:=m_1=m_2=m_3$. Choosing $(p,q,r)=(w_1,w_2,w_3)$ as defined above, we obtain the three rational points $\displaystyle P_i=\left(\frac{1}{v_i^2},w_i\right)$, $i=1,2,3$, in $E'(\Q)$ where $E'$ is defined by $y^2=D x^m+a$.

We consider the automorphisms $\phi_{i}$, $1\le i\le 3$, of $K=\Q(u,v_1,v_2,v_3)$ defined by 
 \begin{eqnarray*}
 \phi_{i}(u)=u, \quad \phi_i(v_i)=-v_i,\quad \phi_i(v_j)=v_j\textrm{ for }i\ne j.
 \end{eqnarray*}
Since $\phi_i(T)=T$ and $\phi_i(D)=D$, the automorphism $\phi_i$ induces a map $\Phi_i$ on $E'$ and its Jacobian $J$.
In particular, $$\Phi_i(P_i)=P_i, \,i=1,2,3,\quad \Phi_i(P_j)=\left(\frac{1}{v_j^2},-w_j\right),\,i\ne j. $$  

We define the $K$-rational divisors $D_i:=(P_i)-(\infty)$, $i=1,2,3$, on $E'$. In accordance with the proof of Theorem \ref{thm1}, the divisors $D_i$ define rational points of infinite order in $J$. In addition, one has $$\Phi_i(D_i)=D_i,\,i=1,2,3,\qquad \Phi_i(D_j)=\left(\left(\frac{1}{v_j^2},-w_j\right)\right)-(\infty)\sim -D_j,\,i\ne j.$$  
We now assume that there are integers $\alpha,\beta,\gamma$ such that $\alpha D_1 +\beta D_2+\gamma D_3\sim 0$. Applying the automorphism $\Phi_1$, one gets $\alpha D_1-\beta D_2-\gamma D_3\sim 0$. Adding the two linear equivalence relations, one obtains $2\alpha D_1\sim 0$, hence $\alpha=0$. Similarly, one applies the automorphism $\Phi_2$ and obtains $\beta=0$. It follows that $D_1$, $D_2$ and $D_3$ are linearly independent. 
 \end{Proof}

\section{A quadratic twist and three twists of even degrees}

Let $f(x)\in\Q[x]$ be a square-free polynomial of degree at least $3$. Fix odd integers $m_1,m_2,m_3\ge 3$. We are investigating the existence of non-square rational numbers $D$ such that the Jacobians of the (hyperelliptic) curves 
\begin{eqnarray*}
    H:Dy^2=f(x),\quad H_1:y^2=x^{m_1}+a D,\quad H_2:y^2=x^{m_2}+b D,\quad H_3:y^2=x^{m_3}+c D
\end{eqnarray*}
are of positive Mordell-Weil rank over $\Q$. 

Setting $M=\operatorname{lcm}(m_1,m_2,m_3)$, and $M_{i}=M/m_i$, $i=1,2,3$. We find the rational number $D$ by solving the following system of equations
\begin{eqnarray}
\label{eq2}
\frac{y_1^2-x_1^{m_1}}{a}=\frac{y_2^2-x_2^{m_2}}{b}=\frac{y_3^2-x_3^{m_3}}{c}=\frac{f(x_4)}{y_4^2},
\end{eqnarray}
where the solutions $x_i,y_i$, $i=1,2,3,4$, satisfy $x_iy_if(x_4)(y_1^2-x_1^{m_1})(y_2^2-x_2^{m_2})(y_3^2-x_3^{m_3})\ne 0$.
To produce a family of parametric solutions to the system (\ref{eq2}),   we set 
\begin{eqnarray*}
x_1&=&v_1^{2}T^{M_1},\, x_2=v_2^2T^{M_2},\, x_3=v_3^2T^{M_3},\, x_4 =u,\\
y_1&=&pT^{(M-1)/2},\, y_2=qT^{(M-1)/2},\,y_3=rT^{(M-1)/2},\,y_4=\frac{1}{T^{(M-1)/2}},\end{eqnarray*} 
where $u,v_1,v_2,v_3$ are rational parameters. After simplification, one obtains
\[T=\frac{p^2-af(u)}{v_1^{2m_1}}=\frac{q^2-bf(u)}{v_2^{2m_2}}=\frac{r^2-cf(u)}{v_3^{2m_3}}.\]
It follows that solving the system (\ref{eq2}) is equivalent to finding rational points on the curve 
\[v_2^{2m_2}v_3^{2m_3}(p^2-af(u))=v_1^{2m_1}v_3^{2m_3}(q^2-bf(u))=v_1^{2m_1}v_2^{2m_2}(r^2-cf(u))\] which is an intersection $\mathcal H_{a,b,c}$ of two quadratic surfaces in $\mathbb{P}^3$ over $\Q(v_1,v_2,v_3)$. In addition, the curve $\mathcal H_{a,b,c}$ possesses the four projective rational points $[p:q:r:s]=[\pm v_1^{m_1}:\pm v_2^{m_2}:\pm v_3^{m_3}:0]$. Thus, $\mathcal H_{a,b,c}$ is an elliptic curve.

\begin{Theorem}
 \label{thm3}
     Let $E$ be an elliptic curve defined by $y^2=f(x)$ where $f \in \Q[x]$ is a polynomial of degree $3$ or $4$. Assume, moreover, that the Mordell-Weil rank of $E$ is positive. Let $m_1,m_2,m_3\ge 3$ be odd integers. Consider the hyperelliptic curves 
     $$E_1:y^2=x^{m_1}+a^2,\qquad E_2: y^2=x^{m_2}+b^2,\qquad E_3:y^2=x^{m_3}+c^2,$$ where $a,b,c$ are non-zero rational numbers.
 Then there exists a rational function $D_{2,2m_1,2m_2,2m_3} \in \Q( v_1,v_2,v_3)$ such that the quadratic twist of the curve $E$ and the Jacobians of the $2m_i$-twists of the curves $E_i$, $i=1,2,3$,
 by $D_{2,2m_1,2m_2,2m_3}$
 have positive Mordell-Weil rank over $\Q(v_1, v_2,v_3)$.
 \end{Theorem}
\begin{Proof}
Following the discussion above, one needs to find a rational point on the elliptic curve
\[\mathcal H_{a^2,b^2,c^2}: v_2^{2m_2}v_3^{2m_3}(p^2-a^2f(u))=v_1^{2m_1}v_3^{2m_3}(q^2-b^2f(u))=v_1^{2m_1}v_2^{2m_2}(r^2-c^2f(u)).\]
Since $E$ is of positive Mordell-Weil rank, we choose $u$ such that $(u,y_u)\in E(\Q)$ is of infinite order. Now, the elliptic curve $\mathcal H_{a^2,b^2,c^2}$ is birational to 
\[x^2-T_{1}^2=y^2-T_{2}^2=z^2-T_{3}^2\] where $T_{1}=av_2^{m_2}v_3^{m_3}y_u$, $T_{2}=bv_1^{m_1}v_3^{m_3}y_u$, $T_{3}=cv_1^{m_1}v_2^{m_2}y_u$. According to Proposition \ref{prop1}, the latter curve is of positive Mordell-Weil rank over $\Q(T_1,T_2,T_3)$. Now, the proof is similar to the proof of Theorem \ref{thm1}.
\end{Proof}
We remark that doubling the point of infinite order $[a y_u:by_u:cy_u:1]\in \mathcal H_{a^2,b^2,c^2}(\Q(v_1,v_2,v_3))$ yields the point $[w_1:w_2:w_3:1]$ given by
{\footnotesize\begin{eqnarray*}
\Big[\frac{a^2 b^2 v_3^{2m_3} y_u + a^2 c^2 v_2^{2m_2} y_u - b^2 c^2 v_1^{2m_1} y_u}{2 a b c v_2^{m_2} v_3^{m_3}} &:&
\frac{a^2 b^2 v_3^{2m_3} y_u - a^2 c^2 v_2^{2m_2}y_u + b^2 c^2 v_1^{2m_1} y_u}{2a b c v_1^{m_1}v_3^{m_3}} \\
&:&
    \frac{-a^2 b^2 v_3^{2m_3} y_u + a^2 c^2 v_2^{2m_2}y_u +
    b^2 c^2 v_1^{2m_1}y_u}{2a b c v_1^{m_1}v_2^{m_2}} : 1
\Big]\in\mathcal H_{a^2,b^2,c^2}(\Q(v_1,v_2,v_3)).
\end{eqnarray*}}

 One obtains the following result. 
\begin{Theorem}
\label{thm4}
 Let $E$ be an elliptic curve defined by $y^2=f(x)$ where $f \in \Q[x]$ is a polynomial of degree $3$ or $4$. Assume, moreover, that the Mordell-Weil rank of $E$ is positive. Let $m\ge 3$ be an odd integer. Consider the hyperelliptic curves 
     $E':y^2=x^{m}+a^2$, where $a$ is a non-zero rational number.
 Then there exists a rational function $D_{2,2m} \in \Q( v_1,v_2,v_3)$ such that the
 quadratic twist of the curve $E$ by $D_{2,2m}$ is of positive Mordell-Weil rank and the Mordell-Weil rank of the Jacobian of the $2m$-twist of the curve $E'$
 by $D_{2,2m}$ over $\Q(v_1, v_2,v_3)$ is at least $3$.
\end{Theorem}
\begin{Proof}
We set $m:=m_1=m_2=m_3$ and $a=b=c$ in Theorem \ref{thm3}. The three $\Q(v_1,v_2,v_3)$-rational points 
$\displaystyle P_i:=\left(v_i^2, w_i T^{(m-1)/2}\right)$, $i=1,2,3$, 
on the curve $y^2=x^m+a D$ are of infinite order. 

We consider the automorphisms $\phi_{i}$, $ i=1,2,3$, of $\Q(v_1,v_2,v_3)$ defined by 
 \begin{eqnarray*}
 \phi_i(v_i)=-v_i,\quad \phi_i(v_j)=v_j\textrm{ for }i\ne j.
 \end{eqnarray*}
 Now the proof follows as in the proof of Theorem \ref{thm2}.
\end{Proof}

\section{One even degree twist and three odd degree twists}

Given the curves \[y^2=x^{m_1}+a,\quad y^2=x^{m_2}+b,\quad y^2=x^{m_3}+c,\quad y^2=x^{m_4}+d,\] 
where $a,b,c,d$ are non-zero rationals and $m_i$, $i=1,2,3,4$, are odd integers, we investigate the existence of a rational function $D$ such that the Jacobians of the following twists are of positive Mordell-Weil rank
\[y^2=Dx^{m_1}+a,\quad y^2=Dx^{m_2}+b,\quad y^2=Dx^{m_3}+c,\quad y^2=x^{m_4}+dD.\]
The rational $D$ can be obtained by solving the system
\[\frac{y_1^2-a}{x_1^{m_1}}=\frac{y_2^2-b}{x_2^{m_2}}=\frac{y_3^2-c}{x_3^{m_3}}=\frac{y_4^2-x_4^{m_4}}{d}\]
where $x_iy_i(y_1^2-a)(y_2^2-b)(y_3^2-c)(y_4^2-x_4^{m_4})\ne0$.
Let $M=\operatorname{lcm}(m_1,m_2,m_3,m_4)$ with $M_i=M/m_i$, $i=1,2,3,4$. We set
\begin{eqnarray*}
x_1&=&\frac{1}{v_1^2T^{M_1}},\quad x_2=\frac{1}{v_2^2T^{M_2}},\quad x_3=\frac{1}{v_3^2T^{M_3}},\quad x_4=v_4^2 T^{m_4}, \\
y_1&=&p,\quad y_2=q,\quad y_3=r,\quad y_4=uT^{(M-1)/2}
\end{eqnarray*}
where $u,v_1,v_2,v_3,v_4$ are rational parameters. One now obtains that 
\[(p^2-a)v_1^{2m_1}T=(q^2-b)v_2^{2m_2}T=(r^2-c)v_3^{2m_3}T=\frac{u^2-v_4^{2m_4}T}{d}.\]
In other words, one needs to solve the system 
\[T=\frac{u^2}{v_4^{2m_4}+d(r^2-c)v_3^{2m_3}},\qquad v_1^{2m_1}(p^2-a)=v_2^{2m_2}(q^2-b)=v_3^{2m_3}(r^2-c).\]
We now obtain the following result.

\begin{Theorem}
\label{thm5}
 Let $m_1,m_2,m_3,m_4\ge 3$ be odd integers. Consider the hyperelliptic curves 
     $$ E_1:y^2=x^{m_1}+a^2,\qquad E_2: y^2=x^{m_2}+b^2,\qquad E_3:y^2=x^{m_3}+c^2,\quad E_4:y^2=x^{m_4}+d$$ where $a,b,c,d$ are non-zero rational numbers.
 Then there exists a rational function $D_{m_1,m_2,m_3,2m_4} \in \Q(u, v_1,v_2,v_3)$ such that the
Jacobian of the $m_i$-twists of the curves $E_i$, $i=1,2,3$, by $D_{m_1,m_2,m_3,2m_4}$ and the Jacobian of the $2m_4$-twist of $E_4$ 
 by $D_{m_1,m_2,m_3,2m_4}$ have positive Mordell-Weil rank over $\Q(u, v_1, v_2,v_3)$.
\end{Theorem}
\begin{Proof}
Proposition \ref{prop1} implies that the elliptic curve $v_1^{2m_1}(p^2-a^2)=v_2^{2m_2}(q^2-b^2)=v_3^{2m_3}(r^2-c^2)$ has positive Mordell-Weil rank over $\Q(u,v_1,v_2,v_3)$. The result now follows as in the proofs of Theorem \ref{thm1} and Theorem \ref{thm3}.
\end{Proof}

\end{document}